\documentclass[reqno]{amsart}

\usepackage{graphicx}
\usepackage{amscd}
\usepackage{amstext}
\usepackage{amsmath}
\usepackage{amsbsy}
\usepackage{amsfonts}
\usepackage{amssymb}

\newtheorem{theorem}{Theorem}[section]
\theoremstyle{plain}

\newtheorem{example}[theorem]{Example}

\newtheorem{lemma}[theorem]{Lemma}

\newtheorem{proposition}[theorem]{Proposition}
\newtheorem{remark}[theorem]{Remark}

\newtheorem*{assumption}{Assumption}
\numberwithin{equation}{section}

\newcommand{\BlackBoxes}{\global\overfullrule5pt}

\BlackBoxes

\renewcommand{\epsilon}{\varepsilon}

\newcommand{\R}{{\mathbb R}}

\begin{document}
\title[Regularity of Finite Dimensional Realizations]{Regularity of finite-dimensional realizations for Evolution Equations}
\author{Damir Filipovi\'c and Josef Teichmann}
\address{Damir Filipovi\'c, Department of Operations Research and Financial Engineering,
Princeton University, Princeton, NJ 08544-5263, USA. Josef
Teichmann, Institute of financial and actuarial mathematics, TU
Vienna, Wiedner Hauptstrasse 8-10, A-1040 Vienna, Austria}
\email{dfilipov@princeton.edu, josef.teichmann@fam.tuwien.ac.at}
\date{November 2, 2001 (first draft); \today\;(this draft)}
\thanks{We thank K.~David Elworthy for bringing this interesting
problem to our attention. The second author thanks the patient
participants of the Vienna Seminar in Differential Geometry, in
particular Andreas Cap, Stefan Haller, Andreas Kriegl and Peter
Michor.}

\begin{abstract} We show that a continuous local semiflow of $ C^k $-maps on
a finite-dimensional $ C^k$-manifold $ M $ with boundary is in
fact a local $C^k$-semiflow on $M$ and can be embedded into a
local $ C^k $-flow around interior points of $M$ under some weak
assumption. This result is applied to an open regularity problem
for finite-dimensional realizations of stochastic interest rate
models.
\end{abstract}
\maketitle

\section{Introduction}

Let $ k \geq 1 $ be given. We consider a Banach space $X$ and a
{\emph{continuous local semiflow $Fl$ of $ C^k $-maps}} on an open
subset $ V \subset X$ , i.e.
\begin{enumerate}
\item There is $ \epsilon >0 $ and $ V \subset X $ open with $ Fl:
[0,\epsilon[ \times V \to X $ a continuous map.
\item $ Fl(0,x) = x $ and $ Fl(s,Fl(t,x)) = Fl(s+t,x) $ for $ s,t,s+t \in
[0,\epsilon[ $ and $ x,Fl(t,x) \in V $.
\item The map $ Fl_t : V \to X $ is $ C^k $ for $ t \in [0,\epsilon[$.
\end{enumerate}
To shorten terminology we say that ``$Fl$'' is a continuous local
semiflow of $C^k$-maps on $X$ if for any $ x \in X $ there is an
open neighborhood $ V \subset X $ of $x$ and a continuous local
semiflow $Fl=Fl^{(V)}$ of $C^k$-maps on $V$, such that
$Fl^{(V_1)}=Fl^{(V_2)}$ on $V_1\cap V_2$. Continuous local
semiflows of $ C^k $-maps appear naturally as mild solutions of
nonlinear evolution equations (see Appendix A). The continuous
local semiflow $Fl$ is called \emph{$ C^k $} or \emph{local
$C^k$-semiflow} if $ Fl:[0,\epsilon[ \times V \to X $ is $ C^k $.

We assume that we are given a finite-dimensional
$C^{k}$-submanifold $M$ with boundary of $ X$ such that $M$ is
{\emph{locally invariant}} for $Fl$, i.e.~for every $ x \in M \cap
V $ there is $ \delta_x \in ]0,\epsilon[ $ such that $ Fl(t,x) \in
M $ for $ 0 \leq t \leq \delta_x $. In this case $ Fl $ restricts
in a small open neighborhood of any point $ x \in M \cap V $ to a
continuous local semiflow of $ C^k $-maps on $ M $, see
Lemma~\ref{lemh} below, where we make the restriction precise. A
continuous local semiflow of $C^k$-maps on $M$ is defined as
above, the same for local $C^k$-semiflows: notice however that the
manifold might have a boundary (for the notions of analysis in
this case see for example \cite{Lan:99}). By $ T_x M $ we denote
the (full) tangent space at $ x \in M $, even at the boundary. By
$ {(T_x M)}_{\geq 0} $ we denote the halfspace of inward pointing
tangent vectors for $ x \in \partial M $. The boundary subspace of
this halfspace is the tangent space of the $ C^k $-submanifold
(without boundary) $ \partial M $, these are the tangent vectors
parallel to the boundary.

We shall prove that the restriction of a continuous local semiflow
of $C^k$-maps $Fl$ to a $C^k$-submanifold with boundary $M$ is
jointly $C^{k}$ and can in particular be embedded in a local
$C^{k}$-flow around any interior point of $M$. We shall apply
classical methods from \cite{MonZip:55} developed to solve the
fifth Hilbert problem. Nevertheless we have to face the difficulty
that $Fl$ is only a continuous {\emph{local semi}}flow. We can
prove the result under a weak assumption, which will always be
satisfied with respect to our applications. This problem arises in
several contexts, for example recently in interest rate theory,
see \cite{fil/tei:01,fil/tei:01b}.

We first cite the classical results from Dean Montgomery and Leo
Zippin \cite{MonZip:55} and draw  a simple conclusion, which
illustrates, what are going to do, namely proving a non-linear
version of Example \ref{exp}.

\begin{theorem}\label{mzthm}
Let $M$ be a finite-dimensional $C^{k}$-manifold and
$Fl:\mathbb{R}\times M\rightarrow M$ a continuous flow of $ C^k $-maps
on $M$, then $ Fl $ is a $ C^k $-flow on $ M $.
\end{theorem}
\begin{example}\label{exp}
Let $S$ be a strongly continuous group on a Banach space $X$ and assume
that $M$ is a locally $S$-invariant finite-dimensional
$C^{k}$-submanifold of $X$. Then $M \subset D(A^{k})$, where $A$
denotes the infinitesimal generator of $S$, and the restriction of $A$
to $M$ is a $C^{k-1}$-vector field on $M$.
\end{example}

The paper is organized as follows. In Section~\ref{secproof} we
prove the extension of Theorem~\ref{mzthm} for continuous local
semiflows of $C^k$-maps. In Section~\ref{secappl} we apply this
result to a problem that arises in connection with stochastic
interest rate models as it has been announced in
\cite{fil/tei:01}: finite-dimensional realizations are highly
regular objects, namely given by submanifolds with boundary of
$D(A^\infty)$, where $A$ is the generator of a strongly continuous
semigroup. The appendix contains a regularity result for the
dependence of solutions to evolution equations on the initial
point.

We end this section by the announced lemma. Let $M$ be a
finite-dimensional $C^k$-submanifold with boundary and let $M$ be
locally invariant for $Fl$, as defined above. We denote by $
{\mathbb{R}}^n_{\geq 0} $ the halfspace $ \{ x \in {\mathbb{R}}^n;
\, x_n \geq 0 \} $, consequently $ {\R}_{\geq 0} $ is the positive
halfline including $ 0 $.
\begin{lemma}\label{lemh}
For every $x\in M\cap V$ there exists an open neighborhood
$V'\subset X$ of $x$ and $\epsilon'>0$, such that $Fl(t,y)\in M$
for all $(t,y)\in [0,\epsilon'[\times (V'\cap M)$.
\end{lemma}
\begin{proof}
Take $ x\in M $ and a $C^k$-submanifold chart $ u: U \subset X \to
X $ with $ u(U \cap M) = \{0\} \times W \subset \{ 0 \} \times
\mathbb{R}^{n}_{\ge 0} $, where $ U\subset \overline{U}\subset V $
is open and $ x \in U $, $ W \subset \mathbb{R}^{n}_{\geq 0} $ is
open, convex. Here $ n $ denotes the dimension of $M$. We may
assume that $u$ has a continuous extension on $\overline{U}$ with
$u(\overline{U\cap M})=u(\overline{U}\cap M)= \{0\} \times
\overline{W}$ by restriction of $ U $. The closure of $ U \cap M $
is taken in $ M $.

For $y\in U$ define the lifetime in $U\cap M$
\[ T(y):=\sup\{ 0<t<\epsilon\mid \;\forall 0\le s<t: \, Fl(s,y)\in U\cap M,\}.\]
By continuity of $Fl$ we have $Fl(T(y),y)\in \overline{U\cap
M}\setminus (U\cap M)$ if $T(y)<\epsilon$. We claim that there exists
an open neighborhood $V'\subset U$ of $x$ in $X$ and $\epsilon'>0$ such
that $T(y)\ge \epsilon'$ for all $y\in V'$. Indeed, otherwise we could
find a sequence $(x_n)$ in $U\cap M$ with $x_n\to x$ and
$\epsilon>T(x_n)\to 0$. But this means that
$u(Fl(T(x_n),x_n))\in\{0\}\times (\overline{W}\setminus W)$ converges
to $u(Fl(0,x))=u(x)\in \{0\} \times W$, a contradiction. Whence the
claim, and the lemma follows.
\end{proof}

\section{The classical proof revisited}\label{secproof}

Since we are treating local questions as differentiability, we can
-- without any restriction -- assume that
$f:[0,\epsilon\lbrack\times V\rightarrow\mathbb{R}^{n}_{\geq 0}$
is a given continuous local semiflow of $ C^k $-maps, where $V$ is
open, convex in $\mathbb{R}^{n}_{\geq 0}$. For the notion of
differentiability on manifolds with boundary see for example
\cite{Lan:99}. We do not make a difference in notation between
right derivatives and derivatives, even though on the boundary
points in space or time, respectively, we only calculate right
derivatives. We shall always assume in this section that $f$ is
continuous and $f(t,.)$ is $C^{k}$ for all
$t\in\lbrack0,\epsilon\lbrack$, for some $k\geq1$. We shall write
$f(t,x)=(f_{1}(t,x),...,f_{n}(t,x))$ for
$(t,x)\in\lbrack0,\epsilon\lbrack\times V$.

\begin{assumption}[crucial]
We assume that for any $x\in V$ there is $\epsilon_x>0$ such that
$D_{x}f(t,x)$ is invertible for $0\leq t\leq\epsilon_{x}$ ($D_xf$
denotes the derivative with respect to $x$).
\end{assumption}
\begin{lemma}
\label{continuity}The mapping $(t,x)\mapsto D_{x}f(t,x)$ is continuous.
\end{lemma}

\begin{proof}
For the proof we proceed from the Baire category Theorem and Lemma
2 of \cite{MonZip:55} on p. 198. We then have the following
result:

Let $Z$ be any compact interval, $V$ the open set in
$\mathbb{R}^{n}_{\geq 0}$ and let $F:Z\times
V\rightarrow\mathbb{R}$ be a continuous real valued function, such
that $F(g,.)$ is $C^{1}$ for any $g\in Z$. Given $a\in V$ and
$1\leq i\leq n$, the set of points $g_{0} \in Z$ such that
$\frac{\partial}{\partial x_{i}}F$ is continuous at $(g_{0},a)$ is
dense in $Z$, even more, the set where it is not continuous is of
first category in $Z$.

Let now $a\in V$ be fixed, then the set of points $t_{0}\in\lbrack
0,\epsilon\lbrack$ such that $f_{ij}:=\frac{\partial}{\partial
x_{i}}f_{j}$ is continuous at $(t_{0},a)$, for all $ 1 \leq i,j
\leq n $, is everywhere dense in $[0,\epsilon\lbrack$. We shall
denote this set by $I_{a}$. In addition the determinant
$\det(f_{ij})$ is continuous at these points, too. We want to show
now that for fixed $a\in V$ the mappings $f_{ij}$ are continuous
at $(0,a)$. Notice that the determinant at any point of continuity
$(t_{0},a)$, with $t_{0}\in I_{a}$ small enough, is bounded away
from zero in a neighborhood.

We fix $a\in V$, then for $ t_0 \in [0,\epsilon[ $
\[
f(t_{0}+h,a+y)=f(t_{0},f_{1}(h,a+y),...,f_{n}(h,a+y))
\]
for $h \geq 0 $ and $y\in\mathbb{R}^{n}_{\geq 0}$, both sufficiently small, hence%
\[
D_xf(t_{0}+h,a+y)=D_xf(t_{0},f(h,a+y))\cdot D_xf(h,a+y).
\]
There is $ t_0 \in I_{a} $ such that $ D_xf(t_{0},z) $ is invertible in
a neighborhood of $ a $, hence
\[
D_xf(t_{0},f(h,a+y))^{-1}\cdot D_xf(t_{0}+h,a+y)=D_xf(h,a+y)
\]
and therefore%
\[
id=\lim_{h\downarrow0,y\rightarrow0}D_xf(t_{0},f(h,a+y))^{-1}\cdot
D_xf(t_{0}+h,a+y)=\lim_{h\downarrow0,y\rightarrow0}D_xf(h,a+y)
\]
by continuity of $D_xf$ at $(t_{0},a)$, continuity of $f$ in both
variables and the continuity of the inversion of matrices. So $ 0 \in
I_a $ for all $ a \in V $.

Now we can conclude for arbitrary $t\in]0,\epsilon\lbrack$ in the
following
way:%
\[
D_xf(t+h,a+y)=D_xf(t,f(h,a+y))\cdot D_xf(h,a+y)
\]
for $h\geq0$ and $y\in\mathbb{R}^{n}_{\geq 0}$ sufficiently small,
hence by continuity
at $(0,a)$%
\[
\lim_{h\downarrow0,y\rightarrow0}D_xf(t+h,a+y)=\lim_{h\downarrow0,y\rightarrow
0}D_xf(t,f(h,a+y))\cdot D_xf(h,a+y)=D_xf(t,a).
\]
For left continuity we apply%
\[
D_xf(t,a+y)=D_xf(h,f(t-h,a+y))\cdot D_xf(t-h,a+y)
\]
for $h\geq0$ and $y\in\mathbb{R}^{n}_{\geq 0}$ sufficiently small,
hence by continuity of $ D_x f $ at $(0,a)$ and $ (0,f(t,a)) $,
the continuity of $D_xf$ in the second variable and the existence
of
the inverse for small $h$%
\[
\lim_{h\downarrow0,y\rightarrow0}D_xf(t-h,a+y)=\lim_{h\downarrow0,y\rightarrow
0}D_xf(h,f(t-h,a+y))^{-1}\cdot D_xf(t,a+y)=D_xf(t,a).
\]
Consequently the desired assertion holds.
\end{proof}

In the next step we shall show that there is a derivative at $0$.

\begin{lemma}\label{time-derivative}
The right-hand derivative $\frac{d}{dt}f(t,x)|_{t=0}$ exists for $x\in
V$, and for small
$h \geq 0$ we have the formula%
\[
f(h,x)-x=\int_{0}^{h}D_x f(t,x)dt\cdot(\frac{d}{dt}f(0,x)).
\]
Moreover, $\frac{d}{dt}f(t,.)|_{t=0}:V\rightarrow\mathbb{R}^{n}$ is
continuous.
\end{lemma}

\begin{proof}
We may differentiate with respect to $x$ under the integral sign by
Lemma \ref{continuity} und
uniform convergence, so%
\begin{align*}
T(h,x) &  :=\int_{0}^{h}f(t,x)dt\\ D_xT(h,x) &
:=\int_{0}^{h}D_xf(t,x)dt.
\end{align*}
By the mean value theorem we obtain%
\[
T(h,y)-T(h,x)=D_xT(h,\widetilde{x})(y-x),
\]
where $\widetilde{x}\in\lbrack x,y \rbrack$. Now we take $y=f(p,x)$, then%
\begin{align*}
T(h,y)-T(h,x) &  =\int_{p}^{h+p}f(t,x)dt-\int_{0}^{h}f(t,x)dt\\ &
=\int_{h}^{h+p}f(t,x)dt-\int_{0}^{p}f(t,x)dt,
\end{align*}
which finally yields%
\[
\frac{1}{p}(\int_{0}^{p}f(t+h,x)dt-\int_{0}^{p}f(t,x)dt)=D_xT(h,\widetilde
{x})[\frac{1}{p}(f(p,x)-x)].
\]
This equation can be solved by joint continuity of $(h,z) \mapsto \frac{1}{h}\int_{0}%
^{h}D_xf(t,z)dt$: we obtain for small $h$ and a compact set in $x$ that
the expression is in a small neighborhood of the identity matrix. So
inversion leads to the desired result and then to the given formula.

The formula asserts again by inversion, that the derivative is
continuous with respect to $x$.
\end{proof}

By the semigroup-property and the chain rule, the result of
Lemma~\ref{time-derivative} can be extended for $0<t<\epsilon$, and the
derivative of $f(\cdot,x)$ exists for all $t\in [0,\epsilon[$. Indeed,
\begin{align*}
\lim_{p\downarrow0}\frac{f(t+p,x)-f(t,x)}{p} &  =\lim_{p\downarrow0}%
\frac{f(p,f(t,x))-f(t,x)}{p}=\frac{d}{dt}f(0,f(t,x))\\
\lim_{p\downarrow0}\frac{f(t,x)-f(t-p,x)}{p} &  =\lim_{p\downarrow0}%
\frac{f(p,f(t-p,x))-f(t-p,x)}{p}\\
&  =\lim_{p\downarrow0}\frac{1}{p}(\int_{0}^{p}D_xf(t,f(t-p,x))dt)\frac{d}%
{dt}f(0,f(t-p,x))\\ &  =\frac{d}{dt}f(0,f(t,x))
\end{align*}
by Lemmas \ref{continuity} and \ref{time-derivative}. Consequently for
small $h\ge 0$
\begin{equation}\label{form-derivative}
\frac{d}{dt}f(t,x)=\frac{d}{dt}f(0,f(t,x))=(\int_{0}^{h}D_xf(s,f(t,x))ds)^{-1}%
\cdot(f(h,x)-x).
\end{equation}
In particular $ (t,x) \mapsto \frac{d}{dt}f(t,x)$ is continuous in both
variables on the whole domain of definition.

\begin{lemma} The semiflow $ f $ is $ C^k $ in both variables.
\end{lemma}

\begin{proof}
If $f(t,.)$ is $C^{k}$ for $t\in\lbrack0,\epsilon\lbrack$, then the $r$-jet%
\[
(f(t,x_{0}),D_{x}f(t,x_{0})\cdot x_{1},...,D_{x}^r f(t,x_{0})\cdot
x_{1}\cdot...\cdot x_{r})
\]
for $(t,x_{0},x_{1},...,x_{r})\in\lbrack0,\epsilon\lbrack\times
V\times\mathbb{R}^{n}\times...\times\mathbb{R}^{n}$ is a local
semiflow of $ C^{k-r} $-maps, for $ 0 \leq r \leq k-1 $. For $
r=1$, the $ 1 $-jet is a continuous, local semiflow of $ C^{k-1}
$-maps, by Lemma \ref{continuity}. Assume that for $ r < k $ the
$r$-jet is a continuous, local semiflow then, by Lemma
\ref{continuity} again, the $(r+1)$-jet is continuous. By
induction
\[
(t,x)\mapsto D_{x}^{r}f(t,x)
\]
is continuous in both variables for $0\leq r\leq k$.

If we apply the above results to the $r$-jet for $r < k$, we conclude
by equation \eqref{form-derivative} that $ D^r_x f(t,x) $ can be $( k -
r)$ times differentiated with respect to the $t$-variable, and these
derivatives are continuous. Hence $f$ is $C^{k}$ in both variables.
\end{proof}

\begin{theorem}\label{ck-theorem}
Let $k\geq1$ be given and let $Fl:[0,\epsilon[ \times U \to M $ be
a local semiflow on a finite-dimensional $C^{k}$-manifold $M$ with
boundary, which satisfies the following conditions:
\begin{enumerate}
\item  The semiflow $Fl:[0,\epsilon[ \times U \to M$ is continuous with
$ U \subset M $ open.

\item  The mapping $Fl(t,.)$ is $C^{k}$.

\item  For fixed $x\in U $ there exists $\epsilon_{x}>0$ such that
$T_{x}Fl(t,.)$ is invertible for $0\leq t\leq\epsilon_{x}$.
\end{enumerate}
Then $ Fl $ is $C^{k}$ and for any $ x_0 \in U \setminus \partial
M$ there is a local $C^k$-flow $ \tilde{Fl}: ]-\delta,\delta[
\times V \to M $ with $ V \subset U \setminus \partial M$ open
around $ x_0 $ and $ \delta \leq \epsilon $ such that $ Fl(y,t) =
\tilde{Fl}(y,t) $ for $ y \in V $ and $ 0 \leq t \leq \delta $.
This also holds for the smooth case ($k=\infty$).
\end{theorem}
\begin{proof}
By the previous lemmas the map $ Fl:[0,\epsilon[ \times U \to M $
is a $ C^k$-semiflow on $M$. We fix $ x_0 \in U \setminus \partial
M$, then there is $ 0 < \delta < \epsilon $ and $ W \subset U $
open in $ M \setminus \partial M $, such that $ (t,x) \mapsto
(t,Fl(t,x)) $ is $C^k$-invertible on $ [0,\delta[ \times W $ by
the $C^k$-inverse function theorem on manifolds with boundary (see
\cite{Lan:99}). We then choose an open neighborhood $ V \subset
\cap_{0 \leq t < \delta} Fl(t,W) $ of $ x_0 $ in $ M \setminus
\partial M $. Therefore we can define $ \tilde{Fl}(-t,y) :=
\{Fl(.,.)\}^{-1}(t,y) $ for $ t \in [0,\delta[ $ and $ y \in V $.
Since this is the unique solution in $ z $ of the $ C^k$-equation
$ Fl(t,z)=y $, we obtain a $ C^k$-map $ \tilde{Fl} $. The flow
property holds by uniqueness, too. Notice that $V$ can be chosen
independent of $k$.
\end{proof}
\begin{remark}
Remark that for evolutions (which correspond in the differentiable case
to time-dependent vector fields) we can pass to the extended phase
space and apply the results thereon.
\end{remark}

\section{Applications}\label{secappl}
The following application has been announced in \cite{fil/tei:01} in
connection with finite-dimensional realizations for stochastic models
of the interest rates. Let $X$ be a Banach space, $S$ a strongly
continuous semigroup on $X$ with infinitesimal generator $A:D(A)\to X$,
and let $P:X\to X$ be a locally Lipschitz map. For $x\in D(A)$ we write
\[ \mu(x):=Ax+P(x).\]
Proposition~\ref{locexiuni} yields the existence of a continuous,
local semiflow $Fl^{\mu}$ of mild solutions to the evolution
equation
\begin{equation}\label{homeq}
\frac{d}{dt}x(t)=\mu(x(t)).
\end{equation}
That is, for every $x_0\in X$ there exists a neighborhood $U$ of $x_0$
in $X$ and $T>0$ such that $Fl^{\mu}\in C([0,T]\times U,X)$ and
\begin{equation}\label{mileq}
Fl^{\mu}(t,x)=S_{t}x+\int_{0}^{t}S_{t-s}P(Fl^{\mu}(s,x))ds,\quad\forall
(t,x)\in [0,T]\times U.
\end{equation}
Now let $k\ge 1$, and $M$ be a finite-dimensional
$C^k$-submanifold with boundary in $X$, which is locally invariant
for $Fl^{\mu}$. Hence, by Lemma~\ref{lemh}, $x_0\in M$ implies
$Fl^{\mu}(t,x)\in M$ for all $(t,x)\in [0,T]\times (U\cap M)$, for
some open neighborhood $U$ of $x_0$ in $X$ and $T>0$. By the
methods of \cite{fil:invariant} (see also
Remark~\ref{fil:addinvariant} below) we obtain that necessarily
$M\subset D(A)$ and
\begin{equation}\label{nag1}
\forall x \in M: \, \mu(x)\in T_xM \text{ and } \forall x \in
\partial M: \, \mu(x) \in {(T_x M)}_{\geq 0},
\end{equation}
since $ \mu $ has to be additionally inward pointing. We now can
strengthen this result.
\begin{theorem}\label{thmCk}
Suppose
\begin{equation}\label{assdiff}
P\in \bigcap_{r=0}^k C^{k-r}(X,D(A^r))
\end{equation}
and $D_x^k P$ is locally Lipschitz continuous. Then $M\subset
D(A^k)$ and $\mu|_M$ is a $C^{k-1}$-vector field on~$M$.
\end{theorem}
\begin{proof}
Let $x_0\in M$, and $U$, $T$ as above. Hence $Fl^{\mu}(t,x)\in M$
for all $(t,x)\in [0,T]\times (U\cap M)$. By the assumptions we
made, Theorem~\ref{thmsmooth} applies and we may assume that
$Fl^{\mu}(t,\cdot)$ is $C^k$ on $U$ for all $t\in [0,T]$. Now let
$x\in U\cap M$. We claim that there exists $\epsilon_x>0$ such
that
\begin{equation}\label{cl1}
\text{$D_xFl^\mu(t,x):T_x M\to T_{Fl^\mu(t,x)}M$ is invertible for
$0\le t\le\epsilon_x$.}
\end{equation}
Indeed, let $y\in X$. The directional derivative $D_xFl^\mu(t,x)y$
is continuous in $t$ on $[0,T]$, see \eqref{psicont}. Hence
\eqref{mileq} and dominated convergence imply that
\[ D_xFl^\mu(t,x)y=S_ty+\int_0^t S_{t-s}D P(Fl^\mu(s,x))
D_xFl^\mu(s,x)y\,ds.\] By the bound \eqref{eq1} we conclude that
\[ \sup_{y\in T_x M, \,\|y\|\le 1} \|D_xFl^\mu(t,x)y-y\|
\le \sup_{y\in T_x M, \,\|y\|\le 1}\|S_t y-y\|+ O(t),\] where
$O(t)\to 0$ for $t\to 0$. Since $T_xM$ is finite-dimensional, it
is easy to see that
\[ \sup_{y\in T_x M, \,\|y\|\le 1}\|S_t y-y\|\to 0\quad\text{for $t\to 0$.}\]
Hence there exists $\epsilon_x>0$ such that $D_xFl^\mu(t,x)$
restricted to $T_x M$ is injective, and hence invertible, for all
$t\in [0,\epsilon_x]$. This yields the claim \eqref{cl1}.

Therefore Theorem~\ref{ck-theorem} applies and $Fl^{\mu}:[0,T]\times
(U\cap M)\to M$ is $C^k$. In particular, since $\mu(x)=\partial_t
Fl^\mu(0,x)$, $\mu|_M$ is a $C^{k-1}$-vector field on $M$, and
$Fl^{\mu}(\cdot,x_0)$ is $C^k$ on $[0,T]$.

>From \eqref{mileq} we have
\[ S_tx_0=Fl^{\mu}(t,x_0)-\int_0^t S_{t-s} P(Fl^{\mu}(s,x_0))\,ds,\quad t\in [0,T]. \]
By \eqref{assdiff}, the integral on the right is $C^k$ in $t\in [0,T]$.
Indeed, we obtain inductively by dominated convergence
\begin{multline*}
 \partial_t^r\int_0^t S_{t-s} P(Fl^{\mu}(s,x_0))\,ds =
 \partial_t^{r-1}P(Fl^\mu(t,x_0))+\partial_t^{r-2}AP(Fl^\mu(t,x_0))\\+
 \cdots+A^{r-1}P(Fl^\mu(t,x_0))+\int_0^t S_{t-s} A^r P(Fl^{\mu}(s,x_0))\,ds,
\end{multline*}
for $r\le k$. We conclude that $S_tx_0$ is $C^k$ in $t\in [0,T]$. But
this means that $x_0\in D(A^k)$ and the theorem is proved.
\end{proof}

We now consider a setup that is given in \cite{fil/tei:01}. Let $W$ be
a connected open set in $X$, $d\ge 1$ and
$\sigma=(\sigma_1,\dots,\sigma_d)$ such that
\begin{description}
\item[(A1)] $P$ and $\sigma_i$ are Banach maps from $X$ into $D(A^\infty)$, for $1\le i\le d$.
\item[(A2)] $\mu,\sigma_1,\dots,\sigma_d$ are pointwise linearly independent on $W\cap D(A^\infty)$.
\end{description}
For the definition of a Banach map see \cite{fil/tei:01,ham:82}. The
Banach map principle (\cite[Theorem~5.6.3]{ham:82}) yields that each
$\sigma_i$ generates a local flow $Fl^{\sigma_i}$ on $X$ with the
following property: for every $x_0\in X$ there exists an open
neighborhood $V$ of $x_0$ in $X$ and $T>0$ such that
\[ Fl^{\sigma_i}\in C^\infty(]-T,T[\times V,X)\quad\text{and}
\quad Fl^{\sigma_i}\in  C^\infty(]-T,T[\times V',D(A^\infty)),\]
where $V':=V\cap D(A^\infty)$ is considered as an open set in
$D(A^\infty)$, and $Fl^{\sigma_i}(\cdot,x)$ is the unique solution of
\[ \frac{d}{dt}x(t)=\sigma_i(x(t)),\quad x(0)=x,\qquad(t,x)\in ]-T,T[\times V.\]
Local invariance for $Fl^{\sigma_i}$ is defined as for $Fl^{\mu}$
above.
\begin{theorem}\label{thmregs}
Let $M\subset W$ be a $(d+1)$-dimensional $C^\infty$-submanifold
with boundary of $X$. If $M$ is locally invariant for $Fl^\mu$,
$Fl^{\sigma_1},\dots,Fl^{\sigma_d}$, then $M$ is a
$C^\infty$-submanifold with boundary of $D(A^\infty)$.
\end{theorem}
\begin{proof}
Theorem~\ref{thmCk} implies that $M\subset D(A^\infty)$ and
$\mu|_M$ is a $C^\infty$-vector field on $M$ with respect to the
given differentiable structure as submanifold with boundary of
$X$. Furthermore $ \sigma_1,\dots,\sigma_d $ restrict to smooth
vector fields on $ M $ and $ \sigma_i(x) \in T_x
\partial M $ for $ x \in
\partial M $, since $\sigma_i(x)$ and $-\sigma_i(x)$ have to be
inward pointing by local invariance. We do also have integral
curves for $ \mu $ and $ \sigma_1,\dots,\sigma_d$ on $
D(A^{\infty})$ which coincide with the respective integral curves
on $X$ on the intersection of the domains of definition if they
start from the same point.

We have to construct submanifold charts for $M$, such that $ M $
is also a submanifold with boundary of $ D(A^{\infty})$. We shall
do this by constructing smooth parametrizations $ \alpha:U \to
D(A^{\infty}) $ for any point $ x_0 \in M $, then we apply
\cite[Lemma 3.1]{fil/tei:01}.

Let $x_0\in M \setminus \partial M$. From \cite[Section
2]{fil/tei:01} we know that the vector field $ \mu $ generates a
smooth local semiflow on $ D(A^{\infty}) $, which coincides
locally by uniqueness of integral curves with the local
$C^{\infty}$-flow $Fl^{\mu_M}$ of $ \mu|_M $ on a neighborhood of
$x_0$. This means in particular that $ t \mapsto Fl^{\mu_M}(t,x_0)
$ is smooth with respect to the topology of $ D(A^{\infty}) $. As
in the proof of \cite[Theorem~3.9]{fil/tei:01} it follows, by
{\bf{(A2)}}, that
\[ \alpha(u,x_0):=Fl^{\sigma_1}_{u_1}\circ\cdots\circ
Fl^{\sigma_d}_{u_d}\circ Fl^{\mu_M}_{u_{d+1}}(x_0):U\to
D(A^{\infty}),\] where $U$ is an open, convex (sufficiently small)
neighborhood of $0$ in $\R^{d+1}$, is a diffeomorphism (it has
maximal rank) to an open submanifold $ N \subset M$ with respect
to the differentiable structure as submanifold with boundary of
$X$. But $ \alpha $ is additionally a smooth parametrization of a
submanifold $ N \subset M \subset D(A^\infty) $, therefore we
constructed for the open subset $ N $ of $ M $ an appropriate
chart as submanifold of $D(A^{\infty})$ by \cite[Lemma
3.1]{fil/tei:01}.

For the boundary points $ x_0 \in \partial M $ the argument is
simpler: first we observe that $ {\sigma}_{i}(x) $ are parallel to
the boundary for $ x \in
\partial M $. In this case it is sufficient to define $ \alpha $ on an open, convex subset
$ U \subset \R^{d} \times \R_{\geq 0} $. Again $ \alpha $ is a
smooth diffeomorphism to an open submanifold with boundary $ N
\subset M $ with respect to the original differentiable structure,
but by \cite[Lemma 3.1]{fil/tei:01} this is also a smooth
parametrization of a $N$ as a submanifold of $D(A^{\infty}) $.

Whence we have constructed submanifold charts with respect to
$D(A^{\infty})$, so $ M \subset D(A^{\infty}) $ is also a
submanifold with boundary of $D(A^{\infty})$.
\end{proof}

\begin{remark}\label{fil:addinvariant}
The Nagumo type consistency results in \cite{fil:invariant} have
been derived for submanifolds without boundary. These results can
be extended to submanifolds with boundary. There are two key
points. First, any submanifold with boundary $M$ can be smoothly
embedded in a submanifold without boundary, say $\tilde{M}$, of
the same dimension. Then the main arguments in
\cite{fil:invariant} carry over: to derive the Nagumo type
consistency conditions at a point $ x \in M $ it is enough to have
local viability of the process with initial point $ x\in M $ in
$M$ (and hence in $\tilde{M}$). Consequently we obtain the Nagumo
type conditions for the whole of $ M $ (including the boundary!).
Second, the coefficients of a diffusion process (i.e.~the
coordinate process) viable in a half space have to satisfy the
appropriate inward pointing conditions at the boundary. We refer
to \cite{Mil:97} for the rigorous analysis.
\end{remark}

\begin{appendix}
\section{Regular Dependence on the Initial Point}\label{regularity}
Let $X$ be a Banach space, $S$ a strongly continuous semigroup on $X$
with infinitesimal generator $A$, and $P:\R_{\ge 0}\times X\to X$ a
continuous map. In this section we shall provide the basic existence,
uniqueness and regularity results for the evolution equation
\begin{equation}\label{inva}
\frac{d}{dt}x(t)=Ax(t)+P(t,x(t)).
\end{equation}
We first recall a classical existence and uniqueness result (see
\cite[Theorem~1.2, Chapter 6]{paz:83}).
\begin{theorem}\label{thmexiuni}
Let $T>0$. Suppose $P:[0,T]\times X\to X$ is uniformly Lipschitz
continuous (with constant $C$) on $X$. Then for every $x\in X$ there
exists a unique mild solution $x(t)$, $t\in [0,T]$, to (\ref{inva})
with $x(0)=x$. If $x(t)$ and $y(t)$ are two mild solutions of
(\ref{inva}) with $x(0)=x$ and $y(0)=y$ then
\begin{equation}\label{masineq}
\sup_{t\in [0,T]}\| x(t)-y(t)\|\le M e^{MCT}\|x-y\|,
\end{equation}
where
\begin{equation}\label{Mdef}
M:=\sup_{t\in [0,T]}\|S_t\|.
\end{equation}
\end{theorem}
There is an immediate local version of Theorem~\ref{thmexiuni}. We say
that $P:\R_{\ge 0}\times X\to X$ is {\emph{locally Lipschitz continuous
on $X$}} if for every $T\ge 0$ and $K\ge 0$ there exists $C=C(T,K)$
such that
\[ \|P(t,x)-P(t,y)\|\le C\|x-y\| \]
for all $t\in [0,T]$, and $x,y\in X$ with $\|x\|\le k$ and $\|y\|\le
k$.
\begin{proposition}\label{locexiuni}
Suppose $P:\R_{\ge 0}\times X\to X$ is locally Lipschitz continuous on
$X$. Let $x_0\in X$. Then there exist a neighborhood $U$ of $x_0$ and
$T>0$ such that, for every $x\in U$, equation (\ref{inva}) has a unique
mild solution $x(t)$, $t\in [0,T]$, with $x(0)=x$. If $x(t)$ and $y(t)$
are two mild solutions of \eqref{inva} with $x(0)=x\in U$ and
$y(0)=y\in U$ then \eqref{masineq} holds, for $M$ as in \eqref{Mdef}
and some $C=C(T,U)$.
\end{proposition}
\begin{proof}
Set $K:=2\|x_0\|$ and fix $T'>0$. Define
\[ \tilde{P}(t,x):=
\begin{cases} P(t,x),&\text{if $\|x\|\le K$,}\\
P(t,Kx/\|x\|),&\text{if $\|x\|> K$.}
\end{cases}\]
Then $\tilde{P}:[0,T']\times X\to X$ is uniformly Lipschitz continuous
on $X$ with constant $C=C(T',K)$. Hence Theorem~\ref{thmexiuni} yields
existence and uniqueness of mild solutions for equation~\eqref{inva}
where $P$ is replace by $\tilde{P}$. By \eqref{masineq} there exists
$0<T\le T'$ and a neighborhood $U$ of $x_0$ such that $\sup_{t\in
[0,T]}\|x(t)\|\le K$ for every mild solution $x(t)$ with $x(0)\in U$.
It is now easy to see that $T$ and $U$ satisfy the assertions of the
proposition.
\end{proof}

Here is the announced regularity result.
\begin{theorem}\label{thmsmooth}
Let $k\ge 1$. Suppose $P:\R_{\ge 0}\times X\to X$ is $C^k$ in $x$,
and $D^k_x P$ is continuous on $\R_{\ge 0}\times X$ and locally
Lipschitz continuous on $X$. Let $x_0\in X$. Then there exists an
open neighborhood $U$ of $x_0$ and $T>0$, and a map $F\in
C([0,T]\times U,H)$ such that, for every $x\in U$, $F(\cdot,x)$ is
the unique mild solution of (\ref{inva}) with $F(0,x)=x$. Moreover
$F(t,\cdot)\in C^k(U,X)$ for all $t\in [0,T]$.
\end{theorem}
\begin{proof}
By assumption, $D_x^r P$ is continuous on $\R_{\ge 0}\times X$ and
locally Lipschitz continuous on $X$, for all $r\le k$. Hence
Proposition~\ref{locexiuni} yields the existence of $U$, $T$ and
$F\in C([0,T]\times U,H)$ such that $F(\cdot,x)\in C([0,T],H)$ is
the unique mild solution of \eqref{inva} with $F(0,x)=x$, for all
$x\in U$. It remains to show regularity of $F(t,\cdot)$.

Let $x\in U$ and $y\in X$. The candidate, say $\psi(t,x,y)$, for the
Gateaux directional derivative $D_x F(t,x)y$ is given by the linear
evolution equation
\begin{equation}\label{lineq}
\begin{aligned}
\frac{d}{dt}\psi(t,x,y)&=A\psi(t,x,y)+D_xP(t,F(t,x))\psi(t,x,y)\\
\psi(0,x,y)&=y.
\end{aligned}
\end{equation}
Since $C_1=C_1(x):=\sup_{t\in [0,T]}\|D_xP(t,F(t,x))\|<\infty$,
Theorem~\ref{thmexiuni} yields the existence of a unique mild
solution
\begin{equation}\label{psicont}
\psi(\cdot,x,y)\in C([0,T],X)
\end{equation}
 to \eqref{lineq}, and by
\eqref{masineq}
\begin{equation}\label{eq1}
\sup_{t\in [0,T]}\|\psi(t,x,y)\|\le M e^{MC_1T}\|y\|.
\end{equation}
Now let $t\in [0,T]$ and $(x_n)$ be a sequence in $U$ converging to
$x$. We claim that
\begin{equation}\label{eqest}
\sup_{y\in X, \,\|y\|\le 1}\|\psi(t,x_n,y)-\psi(t,x,y)\|\to
0,\quad\text{for $n\to\infty$.}
\end{equation}
Indeed, $\Delta_n(t):=\psi(t,x_n,y)-\psi(t,x,y)$ satisfies
\[ \Delta_n(t)=\int_0^t S_{t-s}\left(D_xP(s,F(s,x_n))\psi(s,x_n,y)-D_xP(s,F(s,x))\psi(s,x,y)\right)\,ds.\]
Hence
\[ \|\Delta_n(t)\|\le MC_2\int_0^t \|\Delta_n(s)\|\,ds+M^2C_3e^{M(C_0+C_1)T}\|y\|\|x_n-x\|,\]
where $C_0$ and $C_3$ are local Lipschitz constants of $P$ and $D_xP$,
respectively, and $C_2:=\sup_n\sup_{s\in [0,T]}\|D_xP(s,F(s,x_n))\|$.
By Gronwall's inequality
\[ \|\Delta_n(t)\|\le M^2C_3e^{M(C_0+C_1+C_2)T}\|y\|\|x_n-x\| ,\]
whence \eqref{eqest}.

Next, we claim that
\begin{equation}\label{gdeq}
D_x F(t,x)y=\psi(t,x,y).
\end{equation}
Let $\epsilon_0>0$ be such that $x+\epsilon y\in U$ for all
$\epsilon\in [0,\epsilon_0]$. For such $\epsilon$ we write
$\delta(t,\epsilon):=F(t,x+\epsilon y)-F(t,x)-\epsilon\psi(t,x,y)$, and
obtain
\begin{align*}
\delta(t,\epsilon)&=\int_0^t S_{t-s}(P(s,F(s,x+\epsilon
y))-P(s,F(s,x)))\,ds\\ &\quad-\epsilon \int_0^t
S_{t-s}D_xP(s,F(s,x))\psi(s,x,y)\,ds\\ &=\int_0^t
S_{t-s}\left(D_xP(s,F(s,x))\delta(s,\epsilon)+\Delta(s,\epsilon)\right)\,ds,
\end{align*}
where
\[ \Delta (s,\epsilon):=P(s,F(s,x+\epsilon y))-P(s,F(s,x))-D_xP(s,F(s,x))(F(s,x+\epsilon y)-F(s,x)).\]
By regularity of $P$ and in view of \eqref{masineq} there exists
$C_4=C_4(T,U)$ such that
\[ \sup_{t\in [0,T]}\|\Delta (t,\epsilon)\|\le C_4\epsilon.\]
Hence, writing $C_5:=\sup_{t\in [0,T]}\|D_xP(t,F(t,x))\|$,
\[ \|\delta(t,\epsilon)\|\le C_5MT\int_0^t\|\delta(s,\epsilon)\|\,ds+C_4MT\epsilon, \]
and by Gronwall's inequality $\lim_{\epsilon\to
0}\|\delta(t,\epsilon)\|=0$, whence \eqref{gdeq}.

By \eqref{gdeq} it follows that $D_x F(t,x)y$ is well defined for all
$x\in U$ and $y\in X$, and by \eqref{eqest} the mapping
$D_xF(t,\cdot):U\mapsto L(X)$ is continuous, hence $F(t,\cdot)\in
C^1(U,X)$ for all $t\in [0,T]$.

Higher order regularity is shown by induction of the above argument. We
only sketch the case $C^2$. Let $x\in U$ and $y_1,y_2\in X$, and write
$\psi_2(x,y_1,y_2)$ for the candidate of $D_x^2F(t,x)(y_1,y_2)$, which
solves the inhomogeneous linear evolution equation
\begin{equation}\label{eveq2}
\begin{aligned}
\frac{d}{dt}\psi_2(t,x,y_1,y_2)&=A\psi_2(t,x,y_1,y_2)+D_xP(t,F(t,x))\psi_2(t,x,y_1,y_2)\\
&\quad+D_x^2P(t,F(t,x))(D_x F(t,x)y_1,D_x F(t,x)y_2)\\
\psi_2(0,x,y_1,y_2)&=0.
\end{aligned}
\end{equation}
Notice that the inhomogeneous part, $D_x^2P(t,F(t,x))(D_x F(t,x)y_1,D_x
F(t,x)y_2)$, is continuous in $t\in [0,T]$ by induction. Hence
$\psi_2(\cdot,x,y_1,y_2)\in C([0,T],X)$ is the unique mild solution of
\eqref{eveq2} by Theorem~\ref{thmexiuni}. Now let $t\in [0,T]$. One
shows first that $\psi_2(t,\cdot,y_1,y_2):U\to X$ is continuous,
uniformly in $y_1,y_2\in X$ with $\|y_1\|\le 1$, $\|y_2\|\le 1$ (see
\eqref{eqest}). Then the identity
$D_x^2F(t,x)(y_1,y_2)=\psi_2(t,x,y_1,y_2)$ is proved (see
\eqref{gdeq}), whence $F(t,\cdot)\in C^2(U,X)$.
\end{proof}

\end{appendix}

\end{document}